\documentclass[10pt]{article}
\usepackage{stix2} % stix fonts
 \usepackage[cal=boondoxo]{mathalfa} % mathcal
   \usepackage{amsmath,amsthm,amscd,graphicx}
 \usepackage[pagebackref]{hyperref} 
 \usepackage{color}
\newtheorem{thm}{Theorem}[section]
\newtheorem*{thm*}{Theorem}
\newtheorem*{corr*}{Corollary}
\newtheorem{lemma}[thm]{Lemma}
\newtheorem{prop}[thm]{Proposition}
\newtheorem*{prop*}{Proposition}
\newtheorem{corr}[thm]{Corollary}
\theoremstyle{definition}

\newtheorem*{dfn*}{Definition}

\newtheorem*{exmple*}{Example}

\newtheorem*{conj*}{Conjecture}
\theoremstyle{remark}

\newtheorem{rmk}[thm]{\textit{Remark}}
\newcommand{\bC}{{\mathbb{C}}}

\newcommand{\bP}{{\mathbb{P}}}

\newcommand{\bZ}{{\mathbb{Z}}}
  \newcommand\cO{{\mathcal O}} 
 \newcommand{\comp}{\raise1pt\hbox{{$\scriptscriptstyle\circ$}}}
 \def\lset{\{}  % for { 
\def\rset{\}}  % for } 
\def\set#1{\lset#1\rset} 
\def\st{\mid}   % for | 
\def\sett#1#2{\lset #1 \st #2 \rset}

\def\mapright#1{\mathop{\vbox{\ialign{
                                ##\crcr
    ${\scriptstyle\hfil\;\;#1\;\;\hfil}$\crcr
 \noalign{\kern2pt\nointerlineskip}
    \rightarrowfill\crcr}}\;}}

\def\la{\langle}
\def\ra{\rangle}

\begin{document}

\title{On complex   surfaces with definite intersection form}
\author{Chris Peters}

\maketitle 
\begin{abstract} A compact complex surface  with positive definite intersection lattice is either the projective plane or a false projective plane.
 If the intersection lattice is negative definite, the surface is either a non-minimal secondary Kodaira surface, a non-minimal elliptic surface
  with $b_1=1$,   or a class VII surface with 
 $b_2>0$. In all cases the lattice is odd and diagonalizable over the integers.
  \end{abstract}

\section{Introduction}
The intersection form of a compact connected orientable $4n$-dimensional  manifold is
  bilinear, symmetric and, by Poincar\'e duality, unimodular. 
  As is well known (cf. \cite{milnor58,milnor}) if such a form is indefinite, its isometry class is uniquely determined by the signature
  and type of the form. I recall that the \emph{type} of the form $b$  can be \emph{even} or \emph{odd}. It is even if $b(x,x)$ is even for all elements $x$ of the lattice
  and odd otherwise.  Odd forms are diagonalizable over the integers, but unimodular even forms are evidently not diagonalizable.\footnote{Consider \cite{milnor58,milnor}  for more precise information.}
  
 In the definite situation the situation is dramatically different: the number of isometry classes goes up drastically with the rank.
  See e.g.\ \cite[Ch. V.2.3]{serre}.
  So one might ask whether all definite forms  occur as  intersection forms. This is indeed the case  for  topological manifolds in view of the celebrated result \cite{freedman} by M. Freedman  implying  that every form can be realized as the  intersection
 form of a simply connected compact oriented $4$-manifold. Moreover, its  oriented homeomorphism type is uniquely determined by the intersection 
 form.  
 
 Let me next turn to differentiable $4$-manifolds. For those 
  Donaldson \cite{don83} proved that in the simply connected  situation a definite form
is   diagonalizable. A little later, in \cite{don87bis}  he proved this also for the non-simply connected case.
 Such differentiable $4$-manifolds are easily constructed: take a connected sum of projective planes and projective planes with
 opposite orientation. These  cannot all have a complex structure.  Indeed, e.g.\  \cite[V. Thm. 1.1]{4authors} implies that  the only possible simply connected complex  surface with a positive definite intersection form is the projective plane.
 %However,  this does not cover non-simply connected surfaces. 

 The main results of this note deal with complex surfaces having definite intersection forms. The first of these
 covers the (slightly easier) K\"ahler case. The result reads as follows:
 
 \begin{thm} Let $X$ be a compact  K\"ahler  surface with a definite intersection form, then $X$ is  either the projective plane 
 or a fake projective plane\footnote{Surfaces of the latter sort have been classified in~\cite{fakeplane} and turn out to  have large fundamental groups.}
, that is a surface of general type  with the same Betti numbers as $\bP^2$. In these cases the 
  intersection form  is isometric to the (trivial) odd 
 positive rank $1$ form $(x,y)\mapsto xy$. \label{thm:Main}
 \end{thm}
  
 This result is probably known to experts but I am not aware of any proof in the literature. It leads to a characterisation of fake planes:
 
 \begin{corr} The only  non-simply connected K\"ahler surfaces with a definite intersection form are the fake planes.
 \end{corr}
 
Let me next  consider   non-K\"ahler surfaces. For those the intersection form can also be negative definite. 
In this case, a distinction has to be made   between minimal and non-minimal surfaces, where, I recall, a surface is \emph{minimal} if it does not
 contain exceptional curves, i.e., rational curves of self-intersection $(-1)$.  Non-minimal surfaces are obtained from minimal surfaces by repeatedly blowing up points. Each blowing up introduces an exceptional curve. 
% The only difficult case is that of class VII surfaces with $b_2>0$. These are conjecturally all of Kato-type and then the intersection form must be diagonalizable. As long as  one does not know this, the result
%can be deduced  from  the above (difficult) result of Donaldson. 
The main theorem is as follows:

\begin{thm} Let $X$ be a compact  non-K\"ahler  surface with a definite intersection form. Then either $X$ is a surface of class VII with $b_2>0$,
a non-minimal secondary Kodaira surface, or a blown up properly elliptic surface whose minimal model has invariants $q=b_1=1$ and   $b_2=c_2=0$. In all cases the intersection form is negative definite and diagonalizable (and hence odd).\label{thm:Main2}
\end{thm}
For the (standard) terminology concerning surfaces I refer to \cite[Ch. VI]{4authors}.

\begin{rmk} 1. The  elliptic surfaces in the above theorem have been classified: they are   obtained from the product $\bP_1\times E$, $E$ an elliptic curve, by doing logarithmic transformations in (lifts of) torsion points of $E$ with sum zero.  See \cite[Ch. II.2, Thm. 7.7]{fm2}.
\\
2.  Donaldson's results are not used in  the proof in the K\"ahler case, but instead the Bogomolov--Miyaoka--Yau inequality   (cf.\ \cite[\S VII.4]{4authors})  is invoked.  For the non-K\"ahler situation the Donaldson results  can  likewise be dispensed of provided  the Kato conjecture
holds, i.e.  class VII surfaces have  global spherical shells.
\end{rmk}

\section{Basic facts from surface theory}

 It is well known that the Chern numbers $c_1^2(X)$ and $c_2(X)$ are topological invariants. This is obvious for $c_2$ since it is the Euler number. For $c_1^2$ this is a consequence of a special case of the index theorem \cite[Thm. 8.2.2]{hirz}   which in this case reads
 \begin{equation} \label{eqn:index}
 \tau(X)=\text{index of } X  = \frac 13 (c_1^2(X)   - 2 c_2(X)) .
\end{equation}
Here the index is the index of the intersection form of $X$.
Also Noether's formula (cf.\ \cite[p. 26]{4authors}) is used below. It is a  special case of the Riemann--Roch formula and reads: 
%. In
%the present situation it relates the invariants  $q$ and $p_g$ to topological data   as follows:\index{Noether's formula}
\begin{equation}
\label{eqn:Noether}
 1-q(X) +p_g (X)  
= \frac{1}{12}(c_1^2(X)+c_2(X)),
\end{equation} 
where $q(X)= \dim H^1(X,\cO_X)$ and $p_g= \dim H^2(X,\cO_X)$.
Furthermore, I shall need an expression for the signature of the intersection form in terms of these invariants (cf.\ \cite[Ch. IV.2--3]{4authors}):

\begin{prop} \label{prop:OnIntSect} Let $X$ be a compact complex surface. Then 
\begin{enumerate}
\item $b_1(X)$ is even and equal to $2q(X)$ if and only if $X$ is K\"ahler.
Otherwise $b_1(X)=2q(X)-1$.
\item In the K\"ahler case the signature of the intersection form equals $(2p_g(X)+1, b_2(X)-2p_g(X)-1)$ and  $(2p_g(X), b_2(X)- 2p_g(X)$ otherwise.
\end{enumerate}

\end{prop}

As a consequence,   firstly, $q(X) $ and  $p_g(X)$ are  topological invariants. Secondly,  
for a K\"ahler surface the intersection form $S_X$ can only be indefinite or positive definite while for a non-K\"ahler surface it can 
a priori be 
indefinite, positive definite or negative definite. It is positive definite if and only if $b_2=2p_g\not=0$ and 
negative definite   if  and only if $p_g=0$ and $b_2\not=0$.

The proof of the main results  uses the Enriques--Kodaira classification. I recall  it in the form
in which  it is needed (cf.\ \cite[Ch. VI]{4authors}):
  
  \begin{thm}[Enriques--Kodaira classification]\label{thm:EnrClass}
Every compact  complex   surface belongs to exactly one of the
following classes according to their Kodaira dimension $\kappa$. The invariants $(c_1^2,c_2)$ are given for their
minimal models:
 
\begin{center}
\begin{tabular}{|c | l | c| c |c| c|}
 \hline
$\kappa$ & Class  & &  $q$ &$c_1^2$   &  $c_2$ \\
 \hline
$ -\infty$  &    rational  surfaces &   algebraic  & $ 0$   & $8$ or  $9$  & $4$ or $3$   \\
		 &  ruled surfaces of genus $>0$  & algebraic &  $  g$ & $8(1-g)$&    $4(1-g)$ \\
		  &   class VII surfaces &   non-K\"ahler & $1$  & $-c_2 $ & $\ge       0$\\
		\hline
$0$ &	Two-dimensional tori &   K\"ahler   & $4$ & $0$  & $0$ \\
	         &	K3 surfaces       &   K\"ahler       & $0$            &   $0$     &  $    24$  \\
	        &   primary Kodaira surfaces &  non-K\"ahler      & $2$   &$0$  &  $0$ \\
	        & secondary Kodaira surfaces &non-K\"ahler& 							  $1$    & $0$  &  $0$\\
	         & Enriques surfaces & algebraic & $0$ & $0$   & $12$ \\
	         &bielliptic surfaces &  algebraic & $2$ & $0$  & $0$ \\
	        \hline
$ 1$	  &     properly elliptic surfaces   &   &  & $0$ & $\ge 0$  \\
 \hline
  $  2$  &surfaces of general type  & algebraic &   & $>0$ & $>0$\\
   \hline
\end{tabular}
\end{center} 
 \end{thm}

\section{Proofs of Theorems~\ref{thm:Main} and \ref{thm:Main2}}
 
 Let $X$ be a compact complex surface, $S_X$ the intersection form on
 the free $\bZ$-module  $\mathsf H_X=H^2(X, \bZ)/\text{torsion}$. So $(\mathsf H_X, S_X)$ is the intersection lattice of $X$.

 Let me introduce some further useful notation. The rank $1$ unimodular positive, respectively negative definite lattices are denoted $\la 1\ra$ and $\la -1\ra$
 respectively. The hyperbolic plane $U$ is the rank $2$ lattice with basis $\set{e,f}$ and form (denoted by a dot) 
 given by   $e\cdot e=f\cdot f =0$, $e\cdot f=1$
 For rational  and ruled surfaces the intersection forms are well known: for $\bP^2$ it is $\la 1 \ra$, for the other minimal rational 
 or ruled surfaces it
 is either $\la 1\ra\oplus \la -1\ra$ or $U$. See for example \cite[Prop. II.18, Prop. V.1.]{CAS}.
 So only $\bP^2$ gives a 
 definite intersection form and the other surfaces  can be discarded for the proof of Theorem~\ref{thm:Main}.

 As to minimality, observe the following result:
  \begin{lemma}
  If $X$ is not minimal, then $\mathsf H_X$ is odd and splits off as many  copies of $\la -1\ra$ as blow-ups from a minimal model are needed to obtain $X$. If, moreover $X$ is K\"ahler, $\mathsf H_X$  is indefinite.
 \end{lemma}
  The reason is that if $X$ is not minimal, the class      of an exceptional curve splits off orthogonally whereas a K\"ahler class has positive self-intersection. This makes the latter somewhat easier to handle.
Hence, I first consider the 
\textbf{\emph{K\"ahler case}} where one only has to consider positive definite forms. 
So let me assume this. Then, by   Proposition~\eqref{prop:OnIntSect} one has $\tau=2p_g+1$ . The index theorem~\eqref{eqn:index}
combined with the Noether formula~\eqref{eqn:Noether}
then  yields the following expressions for $c_1^2$ and $c_2$:
\begin{eqnarray*}
c^2_1&= &10p_g-8q+9\\
c_2&= & 2p_g-4q+3 
\end{eqnarray*}
so that $c_1^2-3c_2= 4(p_g+q)$.  The class of surfaces with Kodaira dimension $-\infty$ has already be dealt with. 
From the table of the classification theorem~\ref{thm:EnrClass}, one sees that  for surfaces with Kodaira dimension $0,1$ one has 
$c_1^2-3c_2\le 0$. For surfaces of general type this is the
 Bogomolov--Miyaoka--Yau inequality   (cf.\ \cite[\S VII.4]{4authors}).
 Consequently, $p_q=q=0$ and then 
necessarily $S_X\simeq \la 1 \ra$.

Next, consider the \textbf{\emph{non-K\"ahler surfaces}}. The intersection form can   either be positive definite or negative definite.
  In the former case,  the index equals  $\tau= 2p_g$ and in the latter $\tau= -b_2$ and $p_g=0$.
  From  the list of Theorem~\ref{thm:EnrClass} the surfaces concerned are the class VII surfaces, the Kodaira surfaces and the  properly elliptic surfaces.
\begin{itemize}
\item For minimal class VII surfaces the list shows that   $\tau= \frac 13 (c_1^2- 2c_2) =-c_2\le 0$ and so one only  the negative definite case has to be investigated.
Since  $p_g=0, b_2=c_2$ and so  the intersection form is  negative definite  if and only if $b_2>0$.  Minimal such  surfaces have been constructed by   Inoue in \cite{inoue}.
M.\ Kato has shown in  \cite{ClassVIIshells}  that these admit a    holomorphically embedded  copy of $\sett{z\in \bC^2}{1-\epsilon< |z|<1+\epsilon}$ for some  $\epsilon>0$, and   for which, moreover, the complement in the surface is  connected.
Conversely,   any such  Kato surface, by definition a compact complex surface containing  such  a  so-called ``global spherical shell'' 
must be of class VII and is  a  deformation  of a blown up primary Hopf surface (a complex surface diffeomorphic to $S^3\times S^1$). This implies that the intersection form is  diagonalizable  and
negative definite.  
By   Donaldson's result~\cite{don87bis}  this is true for  any class VII surface with $b_2>0$.\footnote{It is conjectured that 
all class VII surfaces with $b_2>0$ are Kato surfaces, which would prove this directly.}

\item By \cite[Ch V.5]{4authors} minimal Kodaira surfaces either have $b_2=4$ and $p_g=1$ (primary Kodaira surfaces).
These have signature $(2,2)$ and since the form is even, it is isometric to $U\operp U$.  In particular, these  need not be considered,
Else  $b_2=0, p_g=0$ (secondary Kodaira surfaces) with zero intersection form and so  only non-minimal such surfaces have negative definite intersection form.

\item  Minimal non-K\"ahler elliptic surfaces.  Since $c_1^2=0$ and $c_2\ge 0$,  the index theorem~\eqref{eqn:index} shows that $\tau\le  0$ and so we only need to
consider the  the negative definite case. Then $p_g=0$,   and thus $p_g-q+1=- q+1=\frac 1 {12} c_2\ge 0$  implying $q=1$, $b_1=1$,
$c_2=b_2=0$. Again  only non-minimal such surfaces have negative definite diagonalizable intersection form.
\end{itemize}

\begin{rmk} As a  consequence of this result, in the case of compact complex surfaces    the intersection form 
is completely determinable from the  Stiefel--Whitney class class $w_2= c_1\bmod 2$ (this determines whether the form is odd or even),
$c_1$, and the Euler number. So the intersection form does not give supplementary topological information unlike   for  topological manifolds.
It then follows from \cite{freedman} that the  oriented homeomorphism type of a simply connected surface is   uniquely determined by the invariants $w_2,c_1^2$ together with  $c_2$.
It is an open question whether this remain true for any compact complex surface  if one adds the fundamental group to the list  of invariants.
One can at least say that the latter determines whether the surface is K\"ahler or not so that the two classes (K\"ahler or not) can be dealt with separately.
\end{rmk}

\bibliographystyle{acm}

\begin{thebibliography}{10}

\bibitem{4authors}
{\sc Barth, W., Hulek, K., Peters, C., and van~de Ven, A.}
\newblock {\em Compact {C}omplex {S}urfaces}, second enlarged~ed.
\newblock No.~4 in Ergebn. der Math. 3. Folge. Springer-Verlag, Berlin etc.,
  2004.

\bibitem{CAS}
{\sc Beauville, A.}
\newblock {\em Complex algebraic surfaces}, second edition~ed.
\newblock Cambridge University Press, Cambridge, 1996.

\bibitem{don83}
{\sc Donaldson, S.~K.}
\newblock An application of gauge theory to four-dimensional topology.
\newblock {\em J. Differential Geom. 18\/} (1983), 279--315.

\bibitem{don87bis}
{\sc Donaldson, S.~K.}
\newblock The orientation of {Y}ang-{M}ills moduli spaces and {$4$}-manifold
  topology.
\newblock {\em J. Differential Geom. 26}, 3 (1987), 397--428.

\bibitem{freedman}
{\sc Freedman, M.}
\newblock The topology of 4-manifolds.
\newblock {\em J. Diff. Geo. 17\/} (1982), 357--454.

\bibitem{fm2}
{\sc Friedman, R., and Morgan, J.}
\newblock {\em Smooth four-manifolds and complex surfaces}, vol.~27 of {\em
  Ergebnisse der Mathematik und ihrer Grenzgebiete (3) [Results in Mathematics
  and Related Areas (3)]}.
\newblock Springer-Verlag, Berlin, 1994.

\bibitem{hirz}
{\sc Hirzebruch, F.}
\newblock {\em Topological Methods in Algebraic Geometry}, vol.~131 of {\em
  {Grundlehre der Math. Wiss.}}
\newblock Springer Verlag, Berlin etc., 1966.

\bibitem{inoue}
{\sc Inoue, M.}
\newblock New surfaces with no meromorphic functions. {II}.
\newblock In {\em Complex analysis and algebraic geometry}. 1977, pp.~91--106.

\bibitem{ClassVIIshells}
{\sc Kato, M.}
\newblock Compact complex manifolds containing ``global'' spherical shells.
\newblock {\em Proc. Japan Acad. 53}, 1 (1977), 15--16.

\bibitem{milnor58}
{\sc Milnor, J.}
\newblock On simply connected {$4$}-manifolds.
\newblock In {\em Symposium internacional de topolog\'{\i}a algebraica
  {I}nternational symposi um on algebraic topology}. Universidad Nacional
  Aut\'{o}noma de M\'{e}xico and UNESCO, Mexico City, 1958, pp.~122--128.

\bibitem{milnor}
{\sc Milnor, J., and Husemoller, D.}
\newblock {\em {Symmetric Bilinear Forms}}, vol.~73 of {\em Ergebn. der Math.
  und ihrer Grenz\-gebiete. 3. Folge}.
\newblock Springer Verlag, Berlin, 1973.

\bibitem{fakeplane}
{\sc Prasad, G., and Yeung, S.-K.}
\newblock Fake projective planes.
\newblock {\em Invent. Math. 168}, 2 (2007), 321--370.

\bibitem{serre}
{\sc Serre, J.-P.}
\newblock {\em A course in arithmetic}.
\newblock Springer Verlag, Berlin etc., 1973.

\end{thebibliography}

\end{document}